\documentclass[12pt]{article}
\usepackage{amsfonts}
\usepackage{mathrsfs}
\usepackage{amsmath,amssymb}
\baselineskip 5.5pt \lineskip 5.5pt \numberwithin{equation}{section}

\def\ad{\mbox{ad}}

\def\lll{\leftline}
\def\Der{{\rm Der}}
\def\Inn{{\rm Inn}}
\def\Aut{{\rm Aut}}
\def\F{\mathbb{F}}
\def\even{{\bar0}}
\def\odd{{\bar1}}

\def\bo{{\bf1}}
\def\v{\varphi}
\def\Ker{{\rm Ker}}
\def\D{\Delta}
\def\rar{\rightarrow}
\def\UU{\frak{U}}

\def\Sp{{\rm{Span}}}

\def\Der{{\rm Der}}

\def\Im{{\rm Im}}
\def\Inn{{\rm Inn}}
\def\Ker{{\rm Ker}}

\def\rar{\rightarrow}

\def\ma{\mathbb}
\def\vs{\vspace*}
\def\cl{\centerline}

\def\c{\frak{c}}
\def\gg{\frak{g}}

\def\D{\Delta}

\def\LL{\mathcal {L}}
\def\WLL{\widetilde{\mathcal {L}}}
\def\GWLL{\widetilde{\mathcal {L}}[s\Gamma]}

\def\cl{\centerline}

\def\C{\mathbb{C}}

\def\Z{\mathbb{Z}}
\def\QED{\hfill$\Box$\par}
\def\BSQ{\hfill$\blacksquare$}

\newtheorem{theo}{Theorem}[section]
\newtheorem{coro}[theo]{Corollary}
\newtheorem{lemm}[theo]{Lemma}

\begin{document}
\baselineskip 18pt

\cl{\large\bf{The derivation algebra and automorphism group}}
\cl{\large\bf{of the (generalized) twisted $N\!=\!2$ superconformal algebra}
\footnote{Supported by NSF grants 11101056, 11271056 of China}}

\vs{6pt}

\cl{Huanxia Fa$^{1,2)}$}

\cl{\small $^{1)}$Wu Wen-Tsun Key Laboratory of Mathematics and School of Mathematical Sciences,}\vs{-3pt}
\cl{\small University of Science and Technology of China, Hefei 230026, China}

\cl{\small $^{2)}$School of Mathematics and Statistics, Changshu Institute
of Technology, Changshu 215500, China}

\cl{\small E-mail: sd\_huanxia@163.com}

\vs{6pt}

\noindent{\bf{Abstract.}} {In this paper, we determine the
derivation algebra and automorphism group of the twisted $N\!=\!2$ superconformal algebra. Then we generalize the relative results to the generalized twisted $N\!=\!2$ superconformal algebra in the final section.


\noindent{\it{MR(2000) Subject Classification}: 17B05, 17B40, 17B65}\vs{18pt}

\lll{\bf 1\ \ Introduction}
\setcounter{section}{1}\setcounter{theo}{0}

It is well known that the superconformal algebras are closely related the conformal field theory and the string theory and play important roles in both mathematics and physics, which were constructed in \cite{K--AM1997} and \cite{ABS--HLB1976} independently. As for the $N\!=\!2$ superconformal algebras, there are four sectors: the Neveu-Schwarz sector, the Ramond sector, the topological sector and the twisted sector, all of which are closely related to the well-known Virasoro algebra and the super-Virasoro algebra. A series of results have been gained on these algebras
(e.g., \cite{DB--CMP2001,EG--CMP1997,FJS--JMP2007,FLX--AC2011,IK--JFA2004,K--AM1997,
K--LMPA1998,KL--SSWS1998} and the correspondingly cited references).

The {\it twisted $N\!=\!2$ superconformal algebra} $\WLL$ is an infinite-dimensional Lie superalgebra over the complex field $\C$ with the
basis $\{L_n,\,T_r,\,G_p,\,\c\,|\,n\in \Z,\,r\in\frac{1}{2}+\Z,\,p\in\frac{1}{2}\Z\}$, admitting the following non-vanishing super brackets:
\begin{eqnarray}\label{LieB}
\begin{array}{lllll}
&&[L_m,L_n]=(m-n)L_{n+m}+\frac{m^3-m}{12}\delta_{m+n,0}\frak{c},\vs{8pt}\\
&&[L_m,T_r]=-rT_{r+m},\ \ \ \ \ \ \ \ \ \ \ \ \ \ \ \
[T_r, T_s]=\frac{r}{3}\delta_{r+s,0}\frak{c},\vs{8pt}\\
&&[L_m,G_p]=(\frac{m}{2}-p)G_{p+m},\ \ \ \ \ \ \ \ [T_r,G_p]=G_{p+r},\vs{8pt}\\
&&[G_p,G_q]=\left\{\begin{array}{lll}
(-1)^{2p}\big(2L_{p+q}+\frac13(p^2-\frac14)\delta_{p+q,0}\frak{c}\big)
&\mbox{if}\ \,p+q\in\Z,\\[6pt]
(-1)^{2p+1}(p-q)T_{p+q}&\mbox{if}\ \,p+q\in\frac12+\Z.
\end{array}\right.
\end{array}\end{eqnarray}

Obviously, $\WLL$ is ${\mathbb{Z}}_2$-graded:
$\WLL=\WLL_{\overline{0}}\oplus\WLL_{\overline{1}},$ with
\begin{eqnarray*}
\WLL_{\overline{0}}
=\mbox{span}_{\ma C}\{L_m,\,T_r,\,\frak{c}\,|\,m\in\Z,\,r\in\frac{1}{2}+\Z\},
\ \ \ \ \WLL_{\overline{1}}=\mbox{span}_{\ma C}\{G_p\,|\,p\in\frac{1}{2}\Z\},
\end{eqnarray*}
and the Cartan subalgebra  ${\cal H}={\mathbb C}L_0+{\mathbb C}\frak{c}$.
One can easily see that $\WLL$ contains  the well-known Virasoro
algebra $\mathcal {V}ir=$ span$_{\mathbb C} \{L_n,\,\frak{c}\,|\,n\in{\Z}\}$, the super-Virasoro algebras $\mathcal{NS}$ (the $N=1$ Neveu-Schwarz algebra) spanned by
$\{L_n,\,G_{r},\,\frak{c}\,|\,n\in\Z,r\in\frac{1}{2}+\Z\}$, and $\mathcal {R}$
(the $N=1$ Ramond algebra) spanned by $\{L_n,\,G_{n},\,\frak{c}\,|\,n\in\Z\}$.

First we recall some definitions and notations. Let
$\mathfrak{g}=\mathfrak{g}_{\bar{0}}\oplus\mathfrak{g}_{\bar{1}}$ be a Lie superalgebra over $\C$. All elements below are assumed to
be $\Z_2$-homogeneous, where $\Z_2=\{\bar0,\,\bar1\}$. For $x\in\mathfrak{g}$,
we always denote $[x]\in\Z_2$ to be its {\it parity}, i.e.,
$x\in\mathfrak{g}_{[x]}$. A $\Z_2$-homogenous linear map $\mathfrak{d}:\mathfrak{g}\to \mathfrak{g}$ such
that there exists $[\mathfrak{d}]\in \Z_2$, $\mathfrak{d}(\mathfrak{g}_{\bar{i}})\subset \mathfrak{g}_{[\bar{i}+[\mathfrak{d}]]}$ for all $\bar{i}\in\Z_2$ satisfying
\begin{equation}\label{20130110812}
\mathfrak{d}([x,y])=[\mathfrak{d}(x),y]+(-1)^{[\mathfrak{d}][x]}[x,\mathfrak{d}(y)] \ \ \mbox{ \ for \ }x,\,y\in\mathfrak{g},
\end{equation}
is called a Lie superalgebra {\it homogenous derivations of parity
$[\mathfrak{d}]$}. The derivation $\mathfrak{d}$ is called {\it even} if $[\mathfrak{d}]=\bar0$, {\it odd} if $[\mathfrak{d}]=\bar1$. Denote by $\Der_{\bar{i}}(\mathfrak{g})$ the set of homogenous
derivations of parity $\bar{i}$. Then
$\Der(\mathfrak{g})=\Der_{\bar0}(\mathfrak{g})\oplus\Der_{\bar1}(\mathfrak{g})$ is the derivation algebra of $\mathfrak{g}$. Denote by $\ad(\mathfrak{g})$ the inner derivation algebra.

It is easy to see that $\WLL$ is a $\frac12\Z$-graded algebra: $\WLL=\oplus_{p\in\frac12\Z}\WLL_p$ where $\WLL_p=\{x\in\WLL\,|\,[L_0,x]=-p\,x\}$. $\Der(\WLL)$ is also $\frac12\Z$-graded: $\Der(\WLL)=\oplus_{p\in\frac12\Z}\Der_p(\WLL)$ where $\Der_p(\WLL)=\{d\in\Der(\WLL)\,|\,d(\WLL_\alpha)\subseteq\WLL_{\alpha+p},
\,\,\alpha\in\frac12\Z\}$.

The result on the derivation algebra $\Der(\WLL)$ of the twisted $N\!=\!2$ superconformal algebra $\WLL$ can be formulated as the following theorem.

\begin{theo}\label{201302252330}
$\Der(\WLL)={\rm {ad}}(\WLL)$.
\end{theo}

In order to introduce the related corollary, we first present some relative notations.
Firstly, let us recall some related definitions based on a Lie superalgebra $\gg$. Denote by $\tau$ the {\it super-twist map} of $\gg\otimes\gg$: $\tau(x\otimes y)= (-1)^{[x][y]}y\otimes x$ for any $x,\,y\in\gg$ and $\xi$ the {\it super-cyclic map} cyclically permuting the coordinates of $\gg^{\otimes3}$: $\xi=({\bf1}\otimes\tau)\cdot(\tau\otimes\bo):\,x_{1}\otimes x_{2}\otimes x_{3}\mapsto(-1)^{[x_1]([x_2]+[x_3])}x_{2}\otimes x_{3}\otimes x_{1}$ for any $x_i\in \gg,\,i=1,2,3$, where $\bo$ is the identity map of $\gg$. A {\it Lie superalgebra} is a pair $(\gg,\varphi)$ consisting of a vector space $\gg=\gg_\even\oplus\gg_\odd$ and a bilinear map $\varphi:\gg\otimes\gg\to\gg$ satisfying:
\begin{eqnarray*}
&&\varphi(\gg_{\bar{i}},\gg_{\bar{j}})\subset\gg_{\bar{i}+\bar{j}},\ \ \, \Ker(\bo\otimes\bo-\tau)\subset\Ker\,\varphi,\ \ \,
\varphi\cdot(\bo\otimes\varphi)\cdot(\bo\otimes\bo\otimes\bo+\xi+\xi^{2})=0.
\end{eqnarray*}
A {\it Lie super-coalgebra} is a pair $(\gg,\D)$ consisting of a vector space $\gg=\gg_\even\oplus\gg_\odd$ and a linear map $\D:\gg\rar\gg\otimes\gg$
satisfying:
\begin{eqnarray*}
&&\!\!\!\!\!\!\D(\gg_{\bar{i}})\subset\mbox{$\sum\limits_{\bar{j}\in{\Z_2}}$}\gg_{\bar{j}}\otimes
\gg_{\bar{i}-\bar{j}},\ \ \,\Im\,\D\subset\Im(\bo\otimes\bo-\tau),\ \ \,
(\bo\otimes\bo\otimes\bo+\xi +\xi^{2})\cdot(\bo\otimes\D)\cdot\D=0.
\end{eqnarray*}
Now one can give the definition of a Lie super-bialgebra, which is
a triple $(\gg,\v,\D)$ satisfying:
$\mbox{(i)}\ \ (\gg,\v){\mbox{is a Lie superalgebra}}$,
$\mbox{(ii)}\ \ (\gg,\D){\mbox{is a Lie super-coalgebra}}$,
$\mbox{(iii)}\ \ \D\v(x\otimes y)=x\ast\D y-(-1)^{[x][y]}y\ast\D x$,
for any $x,\,y\in\gg$, where
the symbol ``$\ast$'' means the {\it adjoint diagonal action}:
$x\ast(\mbox{$\sum\limits_{i}$}{a_{i} \otimes b_{i}})=
\mbox{$\sum\limits_{i}$}({[x,a_{i}]\otimes b_{i}+
(-1)^{[x][a_i]}a_{i}\otimes[x,b_{i}]}),\ \ \forall\,\,x,\,a_{i},\,b_{i}\in\gg$,
and in general $[x,y]=\v(x\otimes y)$ for $x,y\in\gg$.

Denote by $\UU(\gg)$ the {\it universal enveloping algebra} of
$\gg$. If $r=\mbox{$\sum\limits_{i}$}{a_{i}\otimes b_{i}}\in\gg\otimes\gg$, then the following elements are in $\UU(\gg)\otimes\UU(\gg)\otimes\UU(\gg)$:
$r^{12} =\mbox{$\sum\limits_{i}$}{a_{i}\otimes b_{i}\otimes\bo}=r\otimes\bo$, $r^{23} =\mbox{$\sum\limits_{i}$}{\bo\otimes a_{i}\otimes b_{i}}=\bo\otimes r$,
$r^{13}=\mbox{$\sum\limits_{i}$}{a_{i} \otimes\bo\otimes b_{i}}
=(\bo\otimes\tau)(r\otimes\bo)=(\tau\otimes\bo)(\bo\otimes r)$,
while the following elements are in $\gg\otimes\gg\otimes\gg$:
$[r^{12},r^{23}]=\mbox{$\sum\limits_{i,j}$}a_{i} \otimes[b_{i},a_{j}]\otimes b_{j}$,
$[r^{12},r^{13}]=\mbox{$\sum\limits_{i,j}$}(-1)^{[a_j][b_i]}[a_{i},a_{j}]\otimes b_{i} \otimes b_{j}$,
$[r^{13},r^{23}]=\mbox{$\sum\limits_{i,j}(-1)^{[a_j][b_i]}$}a_{i}\otimes a_{j} \otimes[b_{i},b_{j}]$.

A {\it coboundary super-bialgebra} is a quadruple $(\gg, \v,\D,r),$ where $(\gg,\v,\D)$ is a Lie super-bialgebra and $r\in\Im(\bo\otimes\bo-\tau)\subset\gg\otimes\gg$ such that
$\D=\D_r$ is a {\it coboundary of $r$}, where
$\D_r(x)=(-1)^{[r][x]}x\ast r$ for any $x\in\gg$. Furthermore, a coboundary Lie super-bialgebra $(\gg,\v,\D,r)$ is called {\it triangular} if it satisfies the following {\it classical Yang-Baxter Equation} $c(r):=[r^{12}, r^{13}]+[r^{12}, r^{23}]+[r^{13},r^{23}]=0$.

For any $\gg$-module $V$, denote $V^\gg=\{v\in V\,|\,\gg.v=0\}$.

The following lemma can be found in \cite{LCZ--JGP2012}, which is generalized from the Lie algebra case given in \cite{NT--JPAA2000}.
\begin{lemm}\label{201303010600}
Let $\gg$ be a Lie superalgebra such that $\gg^\gg$, $(\gg\otimes\gg)^\gg$ and $H^1(\gg,\gg)$ are all equal to zero. Then for any one dimensional central extension $\widetilde{\gg}$ of $\gg$, there is a linear embedding from $H^1(\gg,\gg\otimes\gg)$ into $H^1(\widetilde{\gg},\widetilde{\gg}\otimes\widetilde{\gg})$. In particular, if $H^1(\gg,\gg\otimes\gg)=0$, then $H^1(\widetilde{\gg},\widetilde{\gg}\otimes\widetilde{\gg})=0$.
\end{lemm}

Denote $\WLL$ with $\c=0$ by $\LL$. The following theorem is the main result of \cite{FLX--AC2011}.
\begin{theo}\label{201303010601}
$H^1(\LL,\LL\otimes\LL)$ and every Lie super-bialgebra structure on $\LL$ is coboundary triangular.
\end{theo}

Theorem \ref{201302252330} implies $H^1(\WLL,\WLL)=0$ and $H^1(\LL,\LL)=0$. It is not difficult to see that $\LL^\LL=0$ and $(\LL\otimes\LL)^\LL=0$. Combining Lemma \ref{201303010600} and Theorem \ref{201303010601}, we can deduce the following corollary.
\begin{coro}
$H^1(\WLL,\WLL\otimes\WLL)$ and every Lie super-bialgebra structure on $\WLL$ is coboundary triangular.
\end{coro}

Denote by ${\rm Aut}(\WLL)$ and $\Inn(\WLL)$ the automorphism group and inner automorphism group of $\WLL$. For any $\varphi\in{\rm Aut}(\WLL)$ and $x,\,y\in\WLL$, we have
\begin{eqnarray}\label{201302270600}
&&\varphi(\WLL_{\bar{0}})=\WLL_{\bar{0}},\ \ \
\varphi(\WLL_{\bar{1}})=\WLL_{\bar{1}},\ \ \
\varphi([x,y])=[\varphi(x),\varphi(y)].
\end{eqnarray}
It is easy to see that $\Inn(\WLL)$ is generated by
\begin{eqnarray*}
{\rm exp}(l_0\,{\rm ad}\,L_0)\ \ \,{\rm for\ some}\,\,l_0\in\C^*.
\end{eqnarray*}
Then $\Inn(\WLL)\cong\C^*$. For any $n\in\Z_+$, denote $\Z/n\Z$ by $\Z_n$.
 
\begin{theo}\label{201302281330}
${\rm Aut}(\WLL)=\Inn(\WLL)\rtimes\Z_4$.
\end{theo}

The following corollary follows immediately from Theorem \ref{201302281330} and $\Inn(\WLL)\cong\C^*$.
\begin{coro}
${\rm Aut}(\WLL)=\C^*\rtimes\Z_4$.
\end{coro}

\vs{18pt}

\lll{\bf 2\ \ Proof of Theorem \ref{201302252330}}
\setcounter{section}{2}\setcounter{theo}{0}\setcounter{equation}{0}

\noindent{\it Proof of Theorem \ref{201302252330}}\ \ \,It will follow from a series of lemmas.
\begin{lemm}\label{201301110751}
$\Der_{\bar 1}(\WLL)={\rm {ad}}_{\bar 1}(\WLL)$.
\end{lemm}
\noindent{\it Proof}\ \ \,For any ${\mathfrak{d}_p}\in\Der_{\bar 1}(\WLL)\cap\Der_{p}(\WLL)$ with  $p\in\frac12\Z$, we always have ${\mathfrak{d}_p}(\frak{c})=0$.

If $p\in\Z$, we can assume
\begin{eqnarray}
\begin{array}{lll}\label{201302242112}
&&{\mathfrak{d}_p}(L_i)=a_{p,i}G_{p+i},\ \ \ \ {\mathfrak{d}_p}(T_r)=b_{p,r}G_{p+r},\vs{8pt}\\
&&{\mathfrak{d}_p}(G_q)=
\left\{\begin{array}{ll}
c_{p,q}L_{p+q}+e_q\delta_{p+q,0}\frak{c}&\mbox{if}\ \,q\in\Z,\\[8pt]
d_{p,q}T_{p+q}&\mbox{if}\ \,q\in\frac12+\Z,
\end{array}\right.
\end{array}
\end{eqnarray}
where $a_{p,i},b_{p,r},c_{p,q},d_{p,q},e_q\in\C$.

According to the assumption that ${\mathfrak{d}_p}$ is a derivation, we have the following identity:
\begin{eqnarray*}
&&{\mathfrak{d}_p}([T_r,T_s])=[{\mathfrak{d}_p}(T_r),T_s]+[T_r,{\mathfrak{d}_p}(T_s)],
\ \ \ \forall\,\, r,\,s\in\frac12+\Z,
\end{eqnarray*}
which gives
\begin{eqnarray*}
&&b_{p,r}G_{p+r+s}=b_{p,s}G_{p+r+s}.
\end{eqnarray*}
Then $b_{p,r}=b_{p,s}$ for any $r,\,s\in\frac12+\Z$, which implies that $b_{p,r}$ is a constant. For convenience, we denote
\begin{eqnarray}\label{201302242100}
b_{p,r}=b_{p},\ \ \forall\,\,p\in\Z,\,r\in\frac12+\Z.
\end{eqnarray}

According to the identity  ${\mathfrak{d}_p}([L_i,T_r])=[{\mathfrak{d}_p}(L_i),T_r]+[L_i,{\mathfrak{d}_p}(T_r)]$, we obtain
\begin{eqnarray*}
&&rb_{p,r+i}G_{p+r+i}=a_{p,i}G_{p+r+i}+(p+r-\frac i2)b_{p,r}G_{p+r+i}.
\end{eqnarray*}
Comparing the coefficients of $G_{p+r+i}$, one has
\begin{eqnarray*}
&&rb_{p,r+i}=a_{p,i}+(p+r-\frac i2)b_{p,r},
\end{eqnarray*}
which together with \eqref{201302242100}, gives
\begin{eqnarray}\label{201302242108}
a_{p,i}=(\frac i2-p)b_{p},\ \ \forall\,\,p,\,i\in\Z.
\end{eqnarray}
For any $r,\,q\in\frac12+\Z$, we have  ${\mathfrak{d}_p}([T_r,G_q])=[{\mathfrak{d}_p}(T_r),G_q]+[T_r,{\mathfrak{d}_p}(G_q)]$, which gives
\begin{eqnarray*}
&&c_{p,q+r}L_{p+q+r}+e_{q+r}\delta_{p+q+r,0}\frak{c}
-d_{p,q}\frac r3\delta_{p+q+r,0}\frak{c}\\
&&=(-1)^{2(p+r)}b_{p,r}\big(2L_{p+q+r}+\frac{4(p+r)^2-1}{12}\delta_{p+q+r,0}\frak{c}\big).
\end{eqnarray*}
Comparing the coefficients of $L_{p+q+r}$ and $\c$, one has
\begin{eqnarray*}
&&c_{p,q+r}=-2b_{p,r},\\
&&e_{-p}=\frac r3d_{p,-p-r}-\frac{4(p+r)^2-1}{12}b_{p,r},
\end{eqnarray*}
which together with \eqref{201302242100}, give
\begin{eqnarray}
&&c_{p,n}=-2b_{p},\label{201302242101}\\
&&e_{-p}=\frac r3d_{p,-p-r}-\frac{4(p+r)^2-1}{12}b_{p},\label{201302242102}
\end{eqnarray}
for any $p,\,n\in\Z$ and $r\in\frac12+\Z$.

For any $r\in\frac12+\Z$ and $q\in\Z$, the identity ${\mathfrak{d}_p}([T_r,G_q])=[{\mathfrak{d}_p}(T_r),G_q]+[T_r,{\mathfrak{d}_p}(G_q)]$ gives
\begin{eqnarray*}
d_{p,q+r}T_{p+q+r}=(-1)^{2(p+r)+1}(p+r-q)b_{p,r}T_{p+q+r}+rc_{p,q}T_{p+q+r}.
\end{eqnarray*}
Comparing the coefficients of $T_{p+q+r}$ and using \eqref{201302242100} with \eqref{201302242101}, we obtain
\begin{eqnarray*}
d_{p,q+r}=(p-q-r)b_{p},
\end{eqnarray*}
which implies
\begin{eqnarray}\label{201302242106}
&&d_{p,q}=(p-q)b_{p},\ \ \forall\,\,p\in\Z,\,q\in\frac12+\Z.
\end{eqnarray}
Using \eqref{201302242102} and \eqref{201302242106}, one can deduce
\begin{eqnarray}\label{201302242110}
&&e_{-p}=-\frac13(p^2-\frac14)b_{p},\ \ \forall\,\,p\in\Z.
\end{eqnarray}
Combing the identities given in \eqref{201302242100}, \eqref{201302242108}, \eqref{201302242101}, \eqref{201302242106} and \eqref{201302242110}, ${\mathfrak{d}_p}(L_i)$, ${\mathfrak{d}_p}(T_r)$ and ${\mathfrak{d}_p}(G_q)$ referred in \eqref{201302242112} can be respectively rewritten as follows:
\begin{eqnarray}
\begin{array}{lll}\label{201302242116}
&&{\mathfrak{d}_p}(L_i)=(\frac i2-p)b_{p}G_{p+i},\ \ \ {\mathfrak{d}_p}(T_r)=b_{p}G_{p+r},\ \ \ {\mathfrak{d}_p}(\c)=0,\vs{8pt}\\
&&{\mathfrak{d}_p}(G_q)=
\left\{\begin{array}{ll}
-2b_{p}L_{p+q}-\frac13(p^2-\frac14)b_{p}\delta_{p+q,0}\frak{c}
&\mbox{if}\ \,q\in\Z,\\[8pt]
(p-q)b_{p}T_{p+q}&\mbox{if}\ \,q\in\frac12+\Z,
\end{array}\right.
\end{array}
\end{eqnarray}
for any $p,\,i\in\Z$, $r\in\frac12+\Z$ and $q\in\frac12\Z$. Noticing
\begin{eqnarray*}
\begin{array}{lll}
&&{\rm ad}(-b_{p}G_p)(L_i)=(\frac i2-p)b_{p}G_{p+i},\ \ \ \ {\rm ad}(-b_{p}G_p)(T_r)=b_{p}G_{p+r},
\ \ \ {\rm ad}(-b_{p}G_p)(\c)=0,\vs{8pt}\\
&&{\rm ad}(-b_{p}G_p)(G_q)=
\left\{\begin{array}{lll}
-2b_{p}L_{p+q}-\frac13(p^2-\frac14)b_{p}\delta_{p+q,0}\frak{c}
&\mbox{if}\ \,q\in\Z,\\[8pt]
(p-q)b_{p}T_{p+q}&\mbox{if}\ \,q\in\frac12+\Z,
\end{array}\right.
\end{array}
\end{eqnarray*}
for any $p,\,i\in\Z$, $r\in\frac12+\Z$ and $q\in\frac12\Z$, we claim that
\begin{eqnarray}\label{201302251200}
{\mathfrak{d}_p}={\rm ad}(-b_{p}G_p),\ \ \ \forall\,\,p\in\Z.
\end{eqnarray}

If $p\in\frac 12+\Z$, we can assume (for convenience, we still use the same notations)
\begin{eqnarray}\label{201302250606}
\begin{array}{lll}
&&{\mathfrak{d}_p}(L_i)=a_{p,i}G_{p+i},\ \ \ \ \ \ \ \ \ \ {\mathfrak{d}_p}(T_r)=b_{p,r}G_{p+r},\vs{8pt}\\
&&{\mathfrak{d}_p}(G_q)=
\left\{\begin{array}{ll}
c_{p,q}L_{p+q}+e_q\delta_{p+q,0}\frak{c},&\mbox{if}\ q\in\frac12+\Z,\\[8pt]
d_{p,q}T_{p+q},&\mbox{if}\ q\in\Z,
\end{array}\right.
\end{array}
\end{eqnarray}
where
$a_{p,i},\,b_{p,r},\,c_{p,q},\,d_{p,q},\,e_q\in\C$.

Using the following two identities:
\begin{eqnarray*}
&&{\mathfrak{d}_p}([L_i,T_r])=[{\mathfrak{d}_p}(L_i),T_r]+[L_i,{\mathfrak{d}_p}(T_r)],
\ \ \ \forall\,\, i\in\Z,\,r\in\frac12+\Z,\\
&&{\mathfrak{d}_p}([T_r,T_s])=[{\mathfrak{d}_p}(T_r),T_s]+[T_r,{\mathfrak{d}_p}(T_s)],
\ \ \ \forall\,\, r,\,s\in\frac12+\Z,
\end{eqnarray*}
we can deduce (for convenience, we denote $b_{p,\frac 12}$ by $b_{p}$ for any $p\in\frac12+\Z$)
\begin{eqnarray}\label{201302251133}
&&a_{p,i}=(\frac i2-p)b_{p},\ \ \ b_{p,r}=b_{p},
\ \ \ \forall\,\, i\in\Z,\,p,\,r\in\frac12+\Z.
\end{eqnarray}

For any $p,\,q\in\frac12+\Z$, we have ${\mathfrak{d}_p}([L_0,G_q])=[{\mathfrak{d}_p}(L_0),G_q]+[L_0,{\mathfrak{d}_p}(G_q)]$, which gives
\begin{eqnarray*}
&&(p+q)c_{p,q}L_{p+q}-q(c_{p,q}L_{p+q}+e_{q}\delta_{p+q,0}\frak{c})\\
&&=a_{p,0}(-1)^{2p}\big(2L_{p+q}+\frac13(p^2-\frac14)\delta_{p+q,0}\frak{c}\big).
\end{eqnarray*}
Comparing the coefficients of $L_{p+q}$ and $\c$, we have
\begin{eqnarray*}
qc_{p,q}\!\!\!&=&\!\!\!(p+q)c_{p,q}-2pb_{p},\\
pe_{-p}\!\!\!&=&\!\!\!\frac13(\frac14-p^2)a_{p,0},
\end{eqnarray*}
which imply
\begin{eqnarray}\label{201302251136}
&&c_{p,q}=2b_{p},\ \ \ e_{-p}=\frac13(p^2-\frac14)b_{p},
\ \ \ \forall\,\,p,\,q\in\frac12+\Z.
\end{eqnarray}

For any $p\in\frac12+\Z$ and $q\in\Z$, we have ${\mathfrak{d}_p}([L_0,G_q])=[{\mathfrak{d}_p}(L_0),G_q]+[L_0,{\mathfrak{d}_p}(G_q)]$, which gives
\begin{eqnarray*}
qd_{p,q}T_{p+q}=a_{p,0}(-1)^{2p+1}(q-p)T_{p+q}+(p+q)d_{p,q}T_{p+q}.
\end{eqnarray*}
Comparing the coefficients of $T_{p+q}$, we have
\begin{eqnarray*}
qd_{p,q}=a_{p,0}(q-p)+(p+q)d_{p,q},
\end{eqnarray*}
which together with \eqref{201302251133}, forces
\begin{eqnarray}\label{201302251138}
d_{p,q}=(q-p)b_{p},\ \ \ \forall\,\,p\in\frac12+\Z,\,q\in\Z.
\end{eqnarray}

Combing the identities given in \eqref{201302251133}, \eqref{201302251136} and \eqref{201302251138}, ${\mathfrak{d}_p}(L_i)$, ${\mathfrak{d}_p}(T_r)$ and ${\mathfrak{d}_p}(G_q)$ referred in \eqref{201302250606} can be respectively rewritten as follows:
\begin{eqnarray*}
\begin{array}{lll}
&&{\mathfrak{d}_p}(L_i)=(\frac i2-p)b_{p}G_{p+i},\ \ \ {\mathfrak{d}_p}(T_r)=b_{p}G_{p+r},\ \ \ {\mathfrak{d}_p}(\c)=0,\vs{8pt}\\
&&{\mathfrak{d}_p}(G_q)=
\left\{\begin{array}{ll}
2b_{p}L_{p+q}+\frac13(p^2-\frac14)b_{p}\delta_{p+q,0}\frak{c}
&\mbox{if}\ \,q\in\frac12+\Z,\\[8pt]
(q-p)b_{p}T_{p+q}&\mbox{if}\ \,q\in\Z,
\end{array}\right.
\end{array}
\end{eqnarray*}
for any $p,\,i\in\Z$, $r\in\frac12+\Z$ and $q\in\frac12\Z$. Noticing
\begin{eqnarray*}
\begin{array}{lll}
&&{\rm ad}(-b_{p}G_p)(L_i)=(\frac i2-p)b_{p}G_{p+i},\ \ \ \ {\rm ad}(-b_{p}G_p)(T_r)=b_{p}G_{p+r},
\ \ \ {\rm ad}(-b_{p}G_p)(\c)=0,\vs{8pt}\\
&&{\rm ad}(-b_{p}G_p)(G_q)=
\left\{\begin{array}{lll}
2b_{p}L_{p+q}+\frac13(p^2-\frac14)b_{p}\delta_{p+q,0}\frak{c}
&\mbox{if}\ \,q\in\frac12+\Z,\\[8pt]
(q-p)b_{p}T_{p+q}&\mbox{if}\ \,q\in\Z,
\end{array}\right.
\end{array}
\end{eqnarray*}
for any $i\in\Z$, $p,\,r\in\frac12+\Z$ and $q\in\frac12\Z$, we claim that
\begin{eqnarray*}
{\mathfrak{d}_p}={\rm ad}(-b_{p}G_p),\ \ \ \forall\,\,p\in\frac12+\Z,
\end{eqnarray*}
which combining with \eqref{201302251200}, gives
\begin{eqnarray*}
{\mathfrak{d}_p}={\rm ad}(-b_{p}G_p),\ \ \ \forall\,\,p\in\frac12\Z.
\end{eqnarray*}
Then this lemma follows.\QED

\begin{lemm}\label{201302251220}
For any ${\mathfrak{D}_p}\in\Der_{\bar 0}(\WLL)\cap\Der_{p}(\WLL)$ with $p\in\frac12\Z^*$, we always have ${\mathfrak{D}_p}\in{\rm {ad}}(\WLL)$.
\end{lemm}
\noindent{\it Proof}\ \ \,It is clear that ${\mathfrak{D}_p}(\frak{c})=0$ for any ${\mathfrak{D}_p}\in\Der_{\bar 0}(\WLL)\cap\Der_{p}(\WLL)$ and $p\in\frac12\Z^*$.

For any $p,\,i\in\Z$, $r\in\frac12+\Z$ and $q\in\frac12\Z$, we can write
\begin{eqnarray}\label{201302251406}
&&{\mathfrak{D}_p}(L_i)=a_{p,i}L_{p+i}+m_i\delta_{p+i,0}\frak{c},\\
&&{\mathfrak{D}_p}(T_r)=b_{p,r}T_{p+r},
\ \ \ {\mathfrak{D}_p}(G_q)=c_{p,q}G_{p+q},
\end{eqnarray}
where
$a_{p,i},b_{p,r},c_{p,q},m_i\in\C.$
For any $p,\,i\in\Z$, we have the following identity:
\begin{eqnarray*}
&&{\mathfrak{D}_p}([L_i,L_0])=[{\mathfrak{D}_p}(L_i),L_0]+[L_i,{\mathfrak{D}_p}(L_0)],
\end{eqnarray*}
which gives
\begin{eqnarray*}
&&i(a_{p,i}L_{p+i}+m_i\delta_{p+i,0}\frak{c})\\
&&=(p+i)a_{p,i}L_{p+i}+a_{p,0}\big((i-p)L_{p+i}+\frac{i^3-i}{12}\delta_{p+i,0}\frak{c}\big).
\end{eqnarray*}
Comparing the coefficients of $L_{p+i}$ and $\c$, we have
\begin{eqnarray*}
ia_{p,i}\!\!\!&=&\!\!\!(p+i)a_{p,i}+(i-p)a_{p,0},\\
pm_{-p}\!\!\!&=&\!\!\!a_{p,0}\frac{p^3-p}{12},
\end{eqnarray*}
which gives
\begin{eqnarray}\label{201302251400}
a_{p,i}=\frac{p-i}{p}a_{p},\ \ \ \ m_{-p}=\frac{p^2-1}{12}a_{p},
\end{eqnarray}
for any $i\in\Z$ and $p\in\Z^*$, where $a_{p,0}$ is denoted by $a_{p}$ for any $p\in\Z^*$.

According to the fact that ${\mathfrak{D}_p}([L_0,T_r])=[{\mathfrak{D}_p}(L_0),T_r]+[L_0,{\mathfrak{D}_p}(T_r)]$,  we obtain
\begin{eqnarray*}
&&rb_{p,r}T_{p+r}=ra_{p}T_{p+r}+(p+r)b_{p,r}T_{p+r},
\end{eqnarray*}
for any $p\in\Z^*$ and $r\in\frac12+\Z$, which implies
\begin{eqnarray}\label{201302251401}
&&b_{p,r}=-\frac{r}{p}a_{p},\ \ \forall\,\,p\in\Z^*,\,r\in\frac12+\Z.
\end{eqnarray}
The identity ${\mathfrak{D}_p}([L_0,G_q])=[{\mathfrak{D}_p}(L_0),G_q]+[L_0,{\mathfrak{D}_p}(G_q)]$ gives
\begin{eqnarray*}
&&qc_{p,q}G_{p+q}=a_{p}(q-\frac{p}{2})G_{p+q}+(p+q)c_{p,q}G_{p+q},
\end{eqnarray*}
for any $p\in\Z^*$ and $q\in\frac12\Z$, from which we can deduce
\begin{eqnarray}\label{201302251402}
&&c_{p,q}=(\frac12-\frac q p)a_{p},\ \ \forall\,\,p\in\Z^*,\,q\in\frac 12\Z.
\end{eqnarray}

Combing the identities given in \eqref{201302251400}, \eqref{201302251401} and \eqref{201302251402}, ${\mathfrak{D}_p}(L_i)$, ${\mathfrak{D}_p}(T_r)$ and ${\mathfrak{D}_p}(G_q)$ presented in \eqref{201302251406} can be respectively rewritten as follows:
\begin{eqnarray}\label{201302251408}
\begin{array}{lll}
&&{\mathfrak{D}_p}(L_i)=\frac{p-i}{p}a_{p}L_{p+i}
+\frac{p^2-1}{12}a_{p}\delta_{p+i,0}\frak{c},
\ \ \ {\mathfrak{D}_p}(\c)=0,\vs{8pt}\\
&&{\mathfrak{D}_p}(T_r)=-\frac{r}{p}a_{p}T_{p+r},
\ \ \ \ \ \ {\mathfrak{D}_p}(G_q)=(\frac12-\frac q p)a_{p}G_{p+q},
\end{array}
\end{eqnarray}
for any $i\in\Z$, $p\in\Z^*$, $r\in\frac12+\Z$ and $q\in\frac12\Z$. Noticing
\begin{eqnarray*}
\begin{array}{lll}
&&{\rm {ad}}(\frac {a_{p}}{p}L_p)(L_i)=\frac{p-i}{p}a_{p}L_{p+i}
+\frac{p^2-1}{12}a_{p}\delta_{p+i,0}\frak{c},
\ \ \ {\rm {ad}}(\frac {a_{p}}{p}L_p)(\c)=0,\vs{8pt}\\
&&{\rm {ad}}(\frac {a_{p}}{p}L_p)(T_r)=-\frac{r}{p}a_{p}T_{p+r},
\ \ \ \ \ \ {\rm {ad}}(\frac {a_{p}}{p}L_p)(G_q)=(\frac12-\frac q p)a_{p}G_{p+q},
\end{array}
\end{eqnarray*}
for any $i\in\Z$, $p\in\Z^*$, $r\in\frac12+\Z$ and $q\in\frac12\Z$, we claim that
\begin{eqnarray}\label{201302252201}
&&{\mathfrak{D}_p}={\rm {ad}}(\frac {a_{p}}{p}L_p),\ \ \ \forall\,\,p\in\Z^*.
\end{eqnarray}

For any $i\in\Z$, $p,\,r\in\frac12+\Z$ and $q\in\frac12\Z$, we can write
\begin{eqnarray}\label{201302252001}
&&{\mathfrak{D}_p}(L_i)=a_{p,i}T_{p+i},
\ \,{\mathfrak{D}_p}(T_r)=b_{p,r}L_{p+r}+n_r\delta_{p+r,0}\frak{c},
\ \,{\mathfrak{D}_p}(G_q)=c_{p,q}G_{p+q},
\end{eqnarray}
where
$a_{p,i},\,b_{p,r},\,c_{p,q},\,n_r\in\C$. For any $i\in\Z$, we have the following identity:
\begin{eqnarray*}
&&{\mathfrak{D}_p}([L_i,L_0])=[{\mathfrak{D}_p}(L_i),L_0]+[L_i,{\mathfrak{D}_p}(L_0)],
\end{eqnarray*}
which gives
\begin{eqnarray*}
&&ia_{p,i}=a_{p,i}(p+i)-pa_{p,0},\ \ \ \forall\,\,i\in\Z,\,p\in\frac12+\Z,
\end{eqnarray*}
and further implies
\begin{eqnarray}\label{201302251800}
a_{p,i}=a_{p},\ \ \ \forall\,\,i\in\Z,\,p\in\frac12+\Z,
\end{eqnarray}
where $a_{p,0}$ is denoted by $a_{p}$ for any $p\in\frac12+\Z$. The following identity: ${\mathfrak{D}_p}([L_0,T_r])=[{\mathfrak{D}_p}(L_0),T_r]+[L_0,{\mathfrak{D}_p}(T_r)]$ gives
\begin{eqnarray*}
&&r(b_{p,r}L_{p+r}+n_r\delta_{p+r,0}\frak{c})
=b_{p,r}(r+p)L_{p+r}-a_{p}\frac p3\delta_{p+r,0}\frak{c}.
\end{eqnarray*}
Comparing the coefficients of $L_{p+r}$ and $\frak{c}$, we obtain
\begin{eqnarray*}
&&rb_{p,r}=(r+p)b_{p,r},\ \ \ pn_{-p}=\frac p3a_{p},
\end{eqnarray*}
for any $r,\,p\in\frac12+\Z$, which imply
\begin{eqnarray}\label{201302252002}
&&b_{p,r}=0,\ \ \ n_{-p}=\frac13a_{p}.
\end{eqnarray}
For any $p\in\frac12+\Z$ and $q\in\frac12\Z$, we have ${\mathfrak{D}_p}([L_0,G_q])=[{\mathfrak{D}_p}(L_0),G_q]+[L_0,{\mathfrak{D}_p}(G_q)]$, from which we can get
\begin{eqnarray*}
&&qc_{p,q}G_{p+q}=(p+q)c_{p,q}G_{p+q}-a_{p}G_{p+q}.
\end{eqnarray*}
Comparing the coefficients of $G_{p+q}$, we have
\begin{eqnarray*}
&&qc_{p,q}=(p+q)c_{p,q}-a_{p},
\end{eqnarray*}
for any $p\in\frac12+\Z$ and $q\in\frac12\Z$, which implies
\begin{eqnarray}\label{201302252003}
&&c_{p,q}=\frac1p a_{p},\ \ \forall\,\,p\in\frac12+\Z,\,q\in\frac 12\Z.
\end{eqnarray}

Combing the identities given in \eqref{201302251800}, \eqref{201302252002} and \eqref{201302252003}, ${\mathfrak{D}_p}(L_i)$, ${\mathfrak{D}_p}(T_r)$ and ${\mathfrak{D}_p}(G_q)$ referred in \eqref{201302252001} can be respectively rewritten as follows:
\begin{eqnarray*}
\begin{array}{lll}
&&{\mathfrak{D}_p}(L_i)=a_{p}T_{p+i},
\ \ \ {\mathfrak{D}_p}(\c)=0,\vs{8pt}\\
&&{\mathfrak{D}_p}(T_r)=\frac{a_{p}}{3}\delta_{p+r,0}\frak{c},
\ \ \ {\mathfrak{D}_p}(G_q)=\frac{a_{p}}{p}G_{p+q},
\end{array}
\end{eqnarray*}
for any $i\in\Z$, $p,\,r\in\frac12+\Z$ and $q\in\frac12\Z$. Noticing
\begin{eqnarray*}
\begin{array}{lll}
&&{\rm {ad}}(\frac {a_{p}}{p}T_p)(L_i)=a_{p}T_{p+i},
\ \ \ {\rm {ad}}(\frac {a_{p}}{p}T_p)(\c)=0,\vs{8pt}\\
&&{\rm {ad}}(\frac {a_{p}}{p}T_p)(T_r)=\frac{a_{p}}{3}\delta_{p+r,0}\frak{c},
\ \ \ {\rm {ad}}(\frac {a_{p}}{p}T_p)(G_q)=\frac{a_{p}}{p}G_{p+q},
\end{array}
\end{eqnarray*}
for any $i\in\Z$, $p,\,r\in\frac12+\Z$ and $q\in\frac12\Z$, we claim that
\begin{eqnarray}\label{201302252202}
&&{\mathfrak{D}_p}={\rm {ad}}(\frac {a_{p}}{p}T_p),\ \ \ \forall\,\,p\in\frac12+\Z.
\end{eqnarray}
Then this lemma follows from \eqref{201302252201} and \eqref{201302252202}.\QED

\begin{lemm}\label{201302252206}
If ${\mathfrak{D}_0}\in\Der_{\bar 0}(\WLL)\cap\Der_{0}(\WLL)$, then ${\mathfrak{D}_0}\in{\rm {ad}}(\WLL)$.
\end{lemm}
\noindent{\it Proof}\ \ \,For any $i\in\Z$, $r\in\frac12+\Z$ and $q\in\frac12\Z$, we can write
\begin{eqnarray}\label{20130225222}
&&\!\!\!\!\!\!\!\!\!\!\!\!
{\mathfrak{D}_0}(\c)=\alpha_0\frak{c},
\ \,{\mathfrak{D}_0}(L_i)=a_iL_i+m_i\delta_{i,0}\frak{c},
\ \,{\mathfrak{D}_0}(T_r)=b_rT_r,\ \ \ {\mathfrak{D}_0}(G_q)=c_qG_q,
\end{eqnarray}
where $a_i,\,b_r,\,c_q,\,m_i,\,\alpha_0\in\C$. For any $i,\,j\in\Z$, we have the following identity:
\begin{eqnarray*}
&&{\mathfrak{D}_0}([L_i,L_j])=[{\mathfrak{D}_0}(L_i),L_j]+[L_i,{\mathfrak{D}_0}(L_j)],
\end{eqnarray*}
which gives
\begin{eqnarray*}
&&(i-j)(a_{i+j}L_{i+j}+m_{i+j}\delta_{i+j,0}\frak{c})
+\frac{i^3-i}{12}\alpha_0\delta_{i+j,0}\frak{c}\\
&&=(a_i+a_j)\big((i-j)L_{i+j}+\frac{i^3-i}{12}\delta_{i+j,0}\frak{c}\big).
\end{eqnarray*}
Comparing the coefficients of $L_{i+j}$ and $\frak{c}$, we have
\begin{eqnarray*}
&&(i-j)a_{i+j}=(i-j)(a_i+a_j),\\
&&2im_{0}+\frac{i^3-i}{12}{\alpha_0}=(a_i+a_{-i})\frac{i^3-i}{12},
\ \ \ \forall\,\,i,\,j\in\Z,
\end{eqnarray*}
from which we can deduce
\begin{eqnarray}\label{201302252300}
&&m_0={\alpha_0}=0\ \ \ {\rm and}\ \ \ a_i=ia_1,\ \ \ \forall\,\,i\in\Z.
\end{eqnarray}
Then ${\mathfrak{D}_0}(\c)$ and ${\mathfrak{D}_0}(L_i)$ referred in \eqref{20130225222} can be rewritten as
\begin{eqnarray}\label{20130225226}
&&{\mathfrak{D}_0}(\c)=0,
\ \ \ {\mathfrak{D}_0}(L_i)=ia_1L_i,\ \ \ \forall\,\,i\in\Z.
\end{eqnarray}

For any $r,\,s\in\frac12+\Z$, we have the following identity:
\begin{eqnarray*}
&&{\mathfrak{D}_0}([T_r,T_s])=[{\mathfrak{D}_0}(T_r),T_s]+[T_r,{\mathfrak{D}_0}(T_s)], \end{eqnarray*}
which together with \eqref{20130225226}, gives
\begin{eqnarray*}
&&\frac r3(b_r+b_s)\delta_{r+s,0}\frak{c}=0.
\end{eqnarray*}
Then we can deduce
\begin{eqnarray*}
&&b_{-r}=-b_r,\ \ \ \forall\,\,r\in\frac12+\Z.
\end{eqnarray*}
According to the identity ${\mathfrak{D}_0}([L_i,T_r])=[{\mathfrak{D}_0}(L_i),T_r]+[L_i,{\mathfrak{D}_0}(T_r)]$,  we obtain
\begin{eqnarray*}
rb_{i+r}T_{i+r}=r(a_i+b_r)T_{i+r},\ \ \ \forall\,\,i\in\Z,\,r\in\frac12+\Z.
\end{eqnarray*}
Comparing the coefficients of $T_{i+r}$ and recalling \eqref{20130225226}, we have
\begin{eqnarray*}
b_{i+r}=ia_1+b_r,\ \ \ \forall\,\,i\in\Z,\,r\in\frac12+\Z,
\end{eqnarray*}
from which we can deduce
\begin{eqnarray}\label{201302252301}
b_{r}=ra_{1},\ \ \ \forall\,\,r\in\frac12+\Z.
\end{eqnarray}
Then ${\mathfrak{D}_0}(T_r)$ referred in \eqref{20130225222} can be rewritten as
\begin{eqnarray}\label{20130225228}
&&{\mathfrak{D}_0}(T_r)=ra_{1}T_r.
\end{eqnarray}

For any $p,\,q\in\frac12\Z$, we have the following identity:
\begin{eqnarray*}
{\mathfrak{D}_0}([G_{p},G_{q}])
=[{\mathfrak{D}_0}(G_{p}),G_{q}]+[G_{p},{\mathfrak{D}_0}(G_{q})],
\end{eqnarray*}
which gives
\begin{eqnarray*}
\left\{\begin{array}{lll}
&\!\!\!\!\!\!
c_{p}+c_{q}=a_{p+q}&\mbox{if}\ \ \,p+q\in\Z,\\[8pt]
&\!\!\!\!\!\!
c_{p}+c_{q}=b_{p+q}&\mbox{if}\ \ \,p+q\in\frac12+\Z.
\end{array}\right.
\end{eqnarray*}
Then recalling \eqref{201302252300} and \eqref{201302252301}, we can deduce
\begin{eqnarray*}
&&c_q=qa_{1},\ \ \ \forall\,\,q\in\frac12\Z.
\end{eqnarray*}
Then ${\mathfrak{D}_0}(G_q)$ referred in \eqref{20130225222} can be rewritten as
\begin{eqnarray}\label{201302252306}
&&{\mathfrak{D}_0}(G_q)=qa_{1}G_q,\ \ \ \forall\,\,q\in\frac12\Z.
\end{eqnarray}
Then according to \eqref{20130225226}, \eqref{201302252301} and \eqref{201302252306}, we know that ${\mathfrak{D}_0}(T_r)$ and ${\mathfrak{D}_0}(G_q)$ referred in \eqref{20130225222} have been rewritten as
\begin{eqnarray*}
&&{\mathfrak{D}_0}(\c)=0,
\ \ \ {\mathfrak{D}_0}(L_i)=ia_1L_i,\\
&&{\mathfrak{D}_0}(T_r)=ra_{1}T_r,\ \ \ {\mathfrak{D}_0}(G_q)=qa_{1}G_q,
\end{eqnarray*}
for any $i\in\Z$, $r\in\frac12+\Z$ and $q\in\frac12\Z$. Noticing that
\begin{eqnarray*}
&&{\rm {ad}}(-a_{1}L_0)(\c)=0,
\ \ \ {\rm {ad}}(-a_{1}L_0)(L_i)=ia_1L_i,\\
&&{\rm {ad}}(-a_{1}L_0)(T_r)=ra_{1}T_r,\ \ \ {\rm {ad}}(-a_{1}L_0)(G_q)=qa_{1}G_q,
\end{eqnarray*}
for any $i\in\Z$, $r\in\frac12+\Z$ and $q\in\frac12\Z$, we get ${\mathfrak{D}_0}={\rm {ad}}(-a_{1}L_0)$. Then the lemma follows.\QED

By now we have completed the proof of Theorem \ref{201302252330}.\BSQ

\vs{18pt}

\lll{\bf3\ \ Proof of Theorem \ref{201302281330}}
\setcounter{section}{3}\setcounter{theo}{0}\setcounter{equation}{0}

\noindent{\it Proof of Theorem \ref{201302281330}}\ \ \,For any $i\in\Z^*$, $r\in\frac12+\Z$, $q\in\frac12\Z$ and $\sigma\in{\rm Aut}(\WLL)$, we can suppose
\begin{eqnarray}\label{201302270608}
\begin{array}{llll}
&&\sigma(L_0)=a_0L_0+c_0\frak{c},\ \ \ \sigma(\frak{c})=m_0\frak{c},\vs{8pt}\\
&&\sigma(L_i)=\mbox{$\sum\limits_{j\in\Z}$}a_{i,j}L_j+
\mbox{$\sum\limits_{s\in\frac12+\Z}$}b_{i,s}T_s+c_i\frak{c},\vs{8pt}\\
&&\sigma(T_r)=\mbox{$\sum\limits_{j\in\Z}$}d_{r,j}L_j+
\mbox{$\sum\limits_{s\in\frac12+\Z}$}e_{r,s}T_s+f_r\frak{c},\vs{8pt}\\
&&\sigma(G_q)=\mbox{$\sum\limits_{p\in\frac12\Z}$}n_{q,p}G_p,
\end{array}
\end{eqnarray}
where $a_0,\,m_0,\,a_{i,j},\,b_{i,s},\,c_i,\,d_{r,j},\,e_{r,s},\,f_r,\,n_{q,p}\in\C$.

We first claim that $a_0\neq0$. Otherwise, if $a_0=0$, for any $x\in\WLL_p$ with $p\neq0$, we have
\begin{eqnarray*}
-p\sigma(x)=[\sigma(L_0),\sigma(x)]=[c_0\frak{c},\sigma(x)]=0,
\end{eqnarray*}
which is impossible. Thus $a_0\neq0$.

For any $i\in\Z^*$, applying $\sigma$ on $[L_0,L_i]=-iL_i$, we have
\begin{eqnarray*}
&&a_0\big(\mbox{$\sum\limits_{j\in\Z}$}ja_{i,j}L_j
+\mbox{$\sum\limits_{s\in\frac12+\Z}$}sb_{i,s}T_s\big)
=i(\mbox{$\sum\limits_{j\in\Z}$}a_{i,j}L_j
+\mbox{$\sum\limits_{s\in\frac12+\Z}$}b_{i,s}T_s+c_i\frak{c}).
\end{eqnarray*}
Comparing the coefficients of $L_j$, $T_s$ and $\frak{c}$, we obtain
\begin{eqnarray*}
&&(a_0j-i)a_{i,j}=0,\ \ \ (a_0s-i)b_{i,s}=0,\ \ \ ic_i=0,
\end{eqnarray*}
which imply
\begin{eqnarray}\label{201302271200}
&&a_{i,j}=b_{i,s}=c_i=0,\ \ \ \forall\,\,i\in\Z^*,\,j\neq\frac{i}{a_0},\,s\neq\frac{i}{a_0}.
\end{eqnarray}
Furthermore, we claim that $a_0\in\{\pm1,\,\pm2\}$. Otherwise, $\sigma(L_i)$ referred in \eqref{201302270608} can be rewritten as $\sigma(L_i)=0$ for some $i\in\Z^*$, which is impossible. Thus $a_0\in\{\pm1,\,\pm2\}$.

For any $r\in\frac12+\Z$, applying $\sigma$ on $[L_0,T_r]=-rT_r$, we have
\begin{eqnarray*}
&&a_0\big(\mbox{$\sum\limits_{j\in\Z}$}jd_{r,j}L_j
+\mbox{$\sum\limits_{s\in\frac12+\Z}$}se_{r,s}T_s\big)
=r(\mbox{$\sum\limits_{j\in\Z}$}d_{r,j}L_j
+\mbox{$\sum\limits_{s\in\frac12+\Z}$}e_{r,s}T_s+f_r\frak{c}).
\end{eqnarray*}
Comparing the coefficients of $L_j$, $T_s$ and $\frak{c}$, we obtain
\begin{eqnarray*}
&&(a_0j-r)d_{r,j}=0,\ \ \ (a_0s-r)e_{r,s}=0,\ \ \ rf_r=0,
\end{eqnarray*}
which imply
\begin{eqnarray}\label{201302271201}
&&d_{r,j}=e_{r,s}=f_r=0,\ \ \ \forall\,\,r\in\frac12+\Z,\,j\neq\frac{r}{a_0},\,s\neq\frac{r}{a_0}.
\end{eqnarray}
We claim that $a_0\notin\{\pm2\}$. Otherwise, $\sigma(T_r)$ referred in \eqref{201302270608} can be rewritten as $\sigma(T_r)=0$ for any $r\in\frac12+\Z$, which is impossible. Thus $a_0\in\{\pm1\}$.

For any $q\in\frac12\Z$, applying $\sigma$ on $[L_0,G_q]=-qG_q$, we have
\begin{eqnarray*}
&&a_0\mbox{$\sum\limits_{p\in\frac12\Z}$}pn_{q,p}G_p
=q\mbox{$\sum\limits_{p\in\frac12\Z}$}n_{q,p}G_p.
\end{eqnarray*}
Comparing the coefficients of $G_p$, we obtain
\begin{eqnarray*}
&&(a_0p-q)n_{q,p}=0,
\end{eqnarray*}
which together with $a_0\in\{\pm1\}$, implies
\begin{eqnarray}\label{201302271202}
&&n_{q,p}=0,\ \ \ \forall\,\,p\neq\frac{q}{a_0}.
\end{eqnarray}

Combining the identities given in \eqref{201302271200}, \eqref{201302271201} and \eqref{201302271202}, we can rewrite $\sigma(L_i)$, $\sigma(T_r)$ and $\sigma(G_q)$ referred in \eqref{201302270608} as follows:
\begin{eqnarray}\label{201302271210}
&&\sigma(L_i)=a_{i,\frac{i}{a_0}}L_{\frac{i}{a_0}},\ \ \
\sigma(T_r)=e_{r,\frac{r}{a_0}}T_{\frac{r}{a_0}},\ \ \
\sigma(G_q)=n_{q,\frac{q}{a_0}}G_{\frac{q}{a_0}},
\end{eqnarray}
for any $i\in\Z^*$, $r\in\frac12+\Z$, $q\in\frac12\Z$ with $a_0\in\{\pm1\}$.

For the case $a_0=1$, using \eqref{201302271210}, we can write
\begin{eqnarray}\label{201302271909}
\begin{array}{lll}
&&\sigma(\frak{c})=m_0\frak{c},\ \ \
\sigma(L_i)=a_iL_i+c_0\delta_{i,0}\frak{c},\\[8pt]
&&\sigma(T_r)=e_rT_r,\ \ \ \sigma(G_q)=n_qG_q,
\end{array}
\end{eqnarray}
for any $i\in\Z$, $r\in\frac12+\Z$, $q\in\frac12\Z$ and the coefficients are all in $\C$.

For any $i,\,j\in\Z$, applying $\sigma$ to
$[L_i,L_j]=(i-j)L_{i+j}+\frac{i^3-i}{12}\delta_{i+j,0}\frak{c}$, we have
\begin{eqnarray*}
&&a_ia_j(i-j)L_{i+j}+\frac{i^3-i}{12}\delta_{i+j,0}a_ia_j\frak{c}\\
&&=(i-j)a_{i+j}L_{i+j}+(i-j)\delta_{i+j,0}c_0\c+\frac{i^3-i}{12}\delta_{i+j,0}m_0\frak{c}.
\end{eqnarray*}
Comparing the coefficients of $L_{i+j}$ and $\frak{c}$, we obtain
\begin{eqnarray*}
&&(i-j)(a_{i+j}-a_ia_j)=0,\\
&&\frac{i^3-i}{12}a_ia_{-i}=2ic_0+\frac{i^3-i}{12}m_0,
\end{eqnarray*}
for any $i,\,j\in\Z$, from which we can deduce
\begin{eqnarray}\label{201302271900}
&&c_0=0,\ \ \ m_0=1,\ \ \ a_i=a^i_1,\ \ \ \forall\,\,i\in\Z.
\end{eqnarray}
From the identity $\sigma([T_r,T_s])=\frac r3\delta_{r+s,0}\sigma(\c)$, we can deduce
\begin{eqnarray}\label{201302271901}
&&e_re_{-r}=1,\ \ \ \forall\,\,r\in\frac12+\Z.
\end{eqnarray}
Using $\sigma([L_i,T_r])=-r\sigma(T_{r+i})$ and recalling \eqref{201302271900}, we can deduce
\begin{eqnarray}\label{201302271902}
&&e_{r+i}=e_ra^i_1,\ \ \ \forall\,\,i\in\Z,\,r\in\frac12+\Z,
\end{eqnarray}
which together with \eqref{201302271901}, gives
\begin{eqnarray}\label{201302271903}
&&a_{i}=\beta^{2i},\ \ \ e_{r}=\beta^{2r},\ \ \ \forall\,\,i\in\Z,\,r\in\frac12+\Z,
\end{eqnarray}
where $e_{\frac12}$ is denoted by $\beta$ for convenience and also $\beta\neq0$.

By $\sigma([G_{p},G_{q}])=n_{p}n_{q}[G_{p},G_{q}]$, we have
\begin{eqnarray*}
\left\{\begin{array}{lll}
&\!\!\!\!\!\!n_{p}n_{q}=a_{p+q}&\mbox{if}\ \,\,p+q\in\Z,\\[8pt]
&\!\!\!\!\!\!n_{p}n_{q}=e_{p+q}&\mbox{if}\ \,\,p+q\in\frac12+\Z,
\end{array}\right.
\end{eqnarray*}
which together with \eqref{201302271903}, give
\begin{eqnarray*}
&&n_{p}n_{q}=\beta^{2p+2q},\ \ \ \forall\,\,p,\,q\in\frac12\Z.
\end{eqnarray*}
Then we can deduce
\begin{eqnarray}\label{201302271912}
&&n_{q}=\varepsilon\beta^{2q},\ \ \ \forall\,\,q\in\frac12\Z,
\end{eqnarray}
where $\varepsilon^2=1$. Combining \eqref{201302271900}, \eqref{201302271902} and \eqref{201302271912}, we can rewrite \eqref{201302271909} as follows:
\begin{eqnarray}\label{201302272300}
\begin{array}{lll}
&&\sigma(\frak{c})=\frak{c},\ \ \ \sigma(L_i)=\beta^{2i}L_i,\\[8pt]
&&\sigma(T_r)=\beta^{2r}T_r,\ \ \ \sigma(G_q)=\varepsilon\beta^{2q}G_q,
\end{array}
\end{eqnarray}
for any $i\in\Z$, $r\in\frac12+\Z$, $q\in\frac12\Z$ and some $\beta\in\C^*$, $\varepsilon^2=1$. The following identities hold:
\begin{eqnarray}\label{201302272311}
\begin{array}{lll}
&&{\rm exp}\big(-({\rm log}\beta^2){\rm ad}L_0\big)(\frak{c})=\frak{c},
\ \ \ {\rm exp}\big(-({\rm log}\beta^2){\rm ad}L_0\big)(L_i)=\beta^{2i}L_i,\\[8pt]
&&{\rm exp}\big(-({\rm log}\beta^2){\rm ad}L_0\big)(T_r)=\beta^{2r}T_r,
\ \ \ {\rm exp}\big(-({\rm log}\beta^2){\rm ad}L_0\big)(G_q)=\beta^{2q}G_q,
\end{array}
\end{eqnarray}
for any $i\in\Z$, $r\in\frac12+\Z$, $q\in\frac12\Z$.

For convenience, we introduce the following isomorphism:
\begin{eqnarray}\label{201302272302}
\begin{array}{lll}
&&\epsilon(\frak{c})=\frak{c},\ \ \ \epsilon(L_i)=L_i,\\[8pt]
&&\epsilon(T_r)=T_r,\ \ \ \epsilon(G_q)=-G_q,
\end{array}
\end{eqnarray}
for any $i\in\Z$, $r\in\frac12+\Z$, $q\in\frac12\Z$. It is easy to check that $\epsilon^k\in{\rm Aut}(\WLL)/\Inn(\WLL)$ for any $k\in\Z_2$.

The $\sigma$ referred in \eqref{201302272300} can be rewritten as:
\begin{eqnarray}\label{201302281300}
&&\sigma=\epsilon^k\,{\rm exp}\big(-({\rm log}\beta^2){\rm ad}L_0\big),
\end{eqnarray}
for some $\beta\in\C^*$, $k\in\Z_2$ and $\epsilon$ is determined by \eqref{201302272302}.

For the case $a_0=-1$, using \eqref{201302271210}, we can write
\begin{eqnarray}\label{201302272000}
\begin{array}{lll}
&&\sigma(\frak{c})=m_0\frak{c},\ \ \
\sigma(L_i)=a_{i}L_{-i}+c_0\delta_{i,0}\frak{c},\\[8pt]
&&\sigma(T_r)=e_rT_{-r},\ \ \ \sigma(G_q)=n_qG_{-q},
\end{array}
\end{eqnarray}
for any $i\in\Z$, $r\in\frac12+\Z$, $q\in\frac12\Z$ and the coefficients are all in $\C$. We still denote $e_{\frac12}$ as $\beta$. Repeating the corresponding process, \eqref{201302272000} can be simplified as follows:
\begin{eqnarray}\label{201302272306}
\begin{array}{lll}
&&\sigma(\frak{c})=-\frak{c},\ \ \
\sigma(L_i)=-\beta^{2i}L_{-i},\\[8pt]
&&\sigma(T_r)=\beta^{2r}T_{-r},\ \ \ \sigma(G_q)=\omega\beta^{2q}G_{-q},
\end{array}
\end{eqnarray}
for any $i\in\Z$, $r\in\frac12+\Z$, $q\in\frac12\Z$ and some $\beta\in\C^*$, $\omega^2=-1$.

For convenience, we introduce the following two isomorphisms:
\begin{eqnarray}
\begin{array}{lll}\label{201302272318}
&&\varpi(\frak{c})=\frak{-c},\ \ \ \varpi(L_i)=-L_{-i},\\[8pt]
&&\varpi(T_r)=T_{-r},\ \ \ \varpi(G_q)=\sqrt{-1} G_{-q},
\end{array}
\end{eqnarray}
for any $i\in\Z$, $r\in\frac12+\Z$, $q\in\frac12\Z$. It is easy to check that $\varpi^k\in{\rm Aut}(\WLL)/\Inn(\WLL)$ for any $k\in\{4k+1\,|\,\forall\,\,k\in\Z\}\cup\{4k+3\,|\,\forall\,\,k\in\Z\}$.

The $\sigma$ referred in \eqref{201302272000} can be rewritten as:
\begin{eqnarray}\label{201302281303}
&&\sigma=\varpi^k\,{\rm exp}\big(-({\rm log}\beta^2){\rm ad}L_0\big),
\end{eqnarray}
for some $\beta\in\C^*$, $k\in\{4k+1\,|\,\forall\,\,k\in\Z\}\cup\{4k+3\,|\,\forall\,\,k\in\Z\}$ and $\varpi$ is determined by \eqref{201302272318}.

Combining \eqref{201302272300}, \eqref{201302281300}, \eqref{201302272306} and \eqref{201302281303}, we finally arrive at the following conclusion:
\begin{eqnarray}\label{201302281306}
&&\sigma=\varpi^k\,{\rm exp}\big(-({\rm log}\beta^2){\rm ad}L_0\big),
\end{eqnarray}
for some $\beta\in\C^*$, $k\in\Z_4$ and $\varpi$ is determined by \eqref{201302272318}. Then this theorem follows.\BSQ

\vs{18pt}

\lll{\bf 4\ \ The generalized case}
\setcounter{section}{4}\setcounter{theo}{0}\setcounter{equation}{0}

Let $\F$ be a field of characteristic zero with the unit identity element $1$, $\Gamma$ an additive subgroup of $\F$, $0$ the identity element of $\Gamma$, and $s\in\F$ satisfying $s\notin\Gamma$ while $2s\in\Gamma$. Denote $\Gamma^*=\Gamma/\{0\}$, $\Gamma_s=s+\Gamma$ and $s\Gamma=\Gamma\cup\Gamma_s$.

The {\it generalized twisted $N\!=\!2$ superconformal algebra}, denoted by $\GWLL$, is an infinite-dimensional Lie superalgebra over $\F$ with the
basis $\{L_\gamma,\,T_\mu,\,G_u,\,\c\,|\,\gamma\in \Gamma,\,\mu\in \Gamma_s,\,u\in s\Gamma\}$ admitting the following non-vanishing super brackets:
\begin{eqnarray}\label{201303010636}
\begin{array}{lllll}
&&[L_{\gamma_1},L_{\gamma_2}]=({\gamma_1}-{\gamma_2})L_{{\gamma_2}+{\gamma_1}}
+\frac{{\gamma_1}^3-{\gamma_1}}{12}\delta_{{\gamma_2}+{\gamma_1},0}\c,\\[8pt]
&&[L_\gamma,T_\mu\,]=-\mu T_{\mu+\gamma},\ \ \ \ \ \ \ \ \ \ \ \ \ \
[T_{\mu_1},T_{\mu_2}]=\frac{\mu_1}{3}\delta_{{\mu_2}+{\mu_1},0}\c,\vs{8pt}\\
&&[L_\gamma,G_u]=(\frac{\gamma}{2}-u)G_{u+\gamma},\ \ \ \ \ \ \ \ [T_\mu,G_u]=G_{u+\mu},\vs{8pt}\\
&&[G_u,G_v]=\left\{\begin{array}{lllll}
2L_{u+v}+\frac13(u^2-\frac14)\delta_{u+v,0}\frak{c}&\mbox{if}\ \,u,\,v\in\Gamma,\\[8pt]
-2L_{u+v}-\frac13(u^2-\frac14)\delta_{u+v,0}\frak{c}
&\mbox{if}\ \,u,\,v\in\Gamma_s,\\[8pt]
(v-u)T_{u+v}&\mbox{if}\ \,u\in\Gamma,\,v\in\Gamma_s,\\[8pt]
(u-v)T_{u+v}&\mbox{if}\ \,v\in\Gamma,\,u\in\Gamma_s.
\end{array}\right.
\end{array}\end{eqnarray}
It is easy to see that $\GWLL$ is $\Z_2$-graded with $\GWLL=\GWLL_{\bar{0}}\oplus\GWLL_{\bar{1}}$, where
\begin{eqnarray*}
\GWLL_{\bar{0}}\!\!\!&=&\!\!\!\Sp_\F\{L_\gamma,T_\mu,\c\,|\,\gamma\in\Gamma,\,\mu\in \Gamma_s\},\\
\GWLL_{\bar{1}}\!\!\!&=&\!\!\!\Sp_\F\{G_u\,|\,u\in s\Gamma\}.
\end{eqnarray*}

A $\Z_2$-homogenous linear map $\mathfrak{D}:\GWLL\to\GWLL$ such
that there exists $[\mathfrak{D}]\in \Z_2$, $\mathfrak{D}(\GWLL_{\bar{i}})\subset \GWLL_{[\bar{i}+[\mathfrak{D}]]}$ for all $\bar{i}\in\Z_2$ satisfying
\begin{eqnarray*}
\mathfrak{D}([x,y])=[\mathfrak{D}(x),y]+(-1)^{[\mathfrak{D}][x]}[x,\mathfrak{D}(y)] \ \ \mbox{ \ for \ }x,\,y\in\GWLL,
\end{eqnarray*}
is called a Lie superalgebra {\it homogenous derivations of parity
$[\mathfrak{D}]$}. The derivation $\mathfrak{D}$ is called {\it even} if $[\mathfrak{D}]=\bar0$, {\it odd} if $[\mathfrak{D}]=\bar1$. Denote by $\Der_{\bar{i}}(\GWLL)$ the set of homogenous
derivations of parity $\bar{i}$. Then
$\Der(\GWLL)=\Der_{\bar0}(\GWLL)\oplus\Der_{\bar1}(\GWLL)$ is the derivation algebra of $\GWLL$. Denote by $\ad(\GWLL)$ the inner derivation algebra.

It is easy to see that $\GWLL$ is a $s\Gamma$-graded algebra: $\GWLL=\oplus_{u\in s\Gamma}\GWLL_u$ where $\GWLL_u=\{x\in\GWLL\,|\,[L_0,x]=-u\,x\}$. $\Der(\GWLL)$ is also $s\Gamma$-graded: $\Der(\GWLL)=\oplus_{u\in s\Gamma}\Der_u(\GWLL)$ where $$\Der_u(\GWLL)=\{d\in\Der(\GWLL)\,|\,d(\GWLL_v)\subseteq
\GWLL_{v+u},\,v\in s\Gamma\}.$$

The result on the derivation algebra $\Der(\GWLL)$ of the generalized twisted $N\!=\!2$ superconformal algebra $\GWLL$ can be formulated as the following theorem, which is not difficult to be generalized from Theorem \ref{201302252330}.
\begin{theo}\label{201303040600}
$\Der(\GWLL)={\rm {ad}}(\GWLL)\oplus{\rm Hom}_\Z(\Gamma,\F)$.
\end{theo}
\noindent{\it Proof}\ \ \,The following results can be obtained from Lemmas \ref{201301110751} and \ref{201302251220} without essential difference:
\begin{eqnarray*}
&&\Der_{\bar 1}(\GWLL)={\rm {ad}}_{\bar 1}(\GWLL)\ \ \,{\rm and}\ \ \, {\mathfrak{D}_p}\in{\rm {ad}}(\GWLL),
\end{eqnarray*}
for any ${\mathfrak{D}_p}\in\Der_{\bar 0}(\GWLL)\cap\Der_{p}(\GWLL)$ with $0\neq p\in s\Gamma$.

For any $\gamma\in\Gamma$, $\mu\in\Gamma_s$ and $u\in s\Gamma$, we can write
\begin{eqnarray}\label{201303040601}
&&\!\!\!\!\!\!
{\mathfrak{D}_0}(\c)=\alpha_0\frak{c},
\ \ \ {\mathfrak{D}_0}(L_\gamma)=a_\gamma L_\gamma+m_\gamma\delta_{\gamma,0}\frak{c},
\ \ \ {\mathfrak{D}_0}(T_\mu)=b_\mu T_\mu,\ \ \ {\mathfrak{D}_0}(G_u)=c_uG_u,
\end{eqnarray}
where $a_\gamma,\,b_\mu,\,c_u,\,m_\gamma,\,\alpha_0\in\F$.

Nearly repeating the proving process of Lemma \ref{201302252206}, we can obtain the following results:
\begin{eqnarray*}
&&m_0={\alpha_0}=0,\ \ \ a_{\gamma_1+\gamma_2}=a_{\gamma_1}+a_{\gamma_2},\\
&&b_{-{\mu}}=-b_{\mu},\ \ \ b_{{\gamma}+{\mu}}=a_{\gamma}+b_{\mu},\\
&&\left\{\begin{array}{lll}
&\!\!\!\!\!\!
c_{u}+c_{v}=a_{u+v}&\mbox{if}\ \ \,u+v\in\Gamma,\\[8pt]
&\!\!\!\!\!\!
c_{u}+c_{v}=b_{u+v}&\mbox{if}\ \ \,u+v\in\Gamma_s,
\end{array}\right.
\end{eqnarray*}
for all ${\gamma},\,\gamma_1,\,\gamma_2\in\Gamma$, ${\mu}\in\Gamma_s$ and $u,\,v\in s\Gamma$. Furthermore, we can deduce
\begin{eqnarray*}
&&c_{\gamma}=a_{\gamma},\ \ \ b_{\mu}=c_{\mu}=\frac{a_{2\mu}}{2},
\end{eqnarray*}
for all ${\gamma},\,\in\Gamma$, ${\mu}\in\Gamma_s$ and $a_{\gamma_1+\gamma_2}=a_{\gamma_1}+a_{\gamma_2}$ for any $\gamma_1,\,\gamma_2\in\Gamma$.

Then the identities referred in \eqref{201303040601} can be rewritten as
\begin{eqnarray*}
&&\!\!\!\!\!\!
{\mathfrak{D}_0}(L_\gamma)=a_\gamma L_\gamma,
\ \ \ {\mathfrak{D}_0}(T_\mu)=\frac{1}{2}a_{2\mu}T_\mu,\\
&&\!\!\!\!\!\!{\mathfrak{D}_0}(G_\gamma)=a_\gamma G_\gamma,
\ \ \ {\mathfrak{D}_0}(G_\mu)=\frac{1}{2}a_{2\mu}G_\mu,
\ \ \ {\mathfrak{D}_0}(\c)=0,
\end{eqnarray*}
for any $\gamma\in\Gamma$, $\mu\in\Gamma_s$ and $u\in s\Gamma$.

For any $\varphi\in{\rm Hom}_\Z(\Gamma,\F)$, one can define the following derivation $\delta_\varphi$:
\begin{eqnarray*}
&&\!\!\!\!\!\!
\delta_\varphi(L_\gamma)=\varphi(\gamma)L_\gamma,
\ \ \ \delta_\varphi(T_\mu)=\frac{1}{2}\varphi(2\mu)T_\mu,\\ &&\!\!\!\!\!\!\delta_\varphi(G_\gamma)=\varphi(\gamma)G_\gamma,
\ \ \ \delta_\varphi(G_\mu)=\frac{1}{2}\varphi(2\mu)G_\mu,
\ \ \ \delta_\varphi(\c)=0,
\end{eqnarray*}
for any $\gamma\in\Gamma$ and $\mu\in\Gamma_s$.

By now we have completed the proof of Theorem \ref{201303040600}.\BSQ

Denote by ${\rm Aut}(\GWLL)$ and $\Inn(\GWLL)$ the automorphism group and inner automorphism group of $\GWLL$. For any $\varphi\in{\rm Aut}(\GWLL)$ and $x,\,y\in\GWLL$, we have
\begin{eqnarray}\label{201302270600}
&&\varphi(\GWLL_{\bar{0}})=\GWLL_{\bar{0}},\ \ \
\varphi(\GWLL_{\bar{1}})=\GWLL_{\bar{1}},\ \ \
\varphi([x,y])=[\varphi(x),\varphi(y)].
\end{eqnarray}
It is easy to see that $\Inn(\GWLL)$ is generated by
\begin{eqnarray*}
{\rm exp}(l_0\,{\rm ad}\,L_0)\ \ \,{\rm for\ some}\,\,l_0\in\F^*.
\end{eqnarray*}
Then $\Inn(\GWLL)\cong\F^*$.

Denote by $\frak{G}$ the subgroup of $\Aut(\GWLL)$, which is generated by the automorphisms determined by \eqref{201303050900} and \eqref{201303050901}. Then the result on the automorphism group $\Aut(\GWLL)$ of the generalized twisted $N\!=\!2$ superconformal algebra $\GWLL$ can be formulated as the following theorem.

\begin{theo}\label{201303050500}
${\rm Aut}(\GWLL)=\Inn(\GWLL)\rtimes\frak{G}$.
\end{theo}

\noindent{\it Proof}}\ \ \,For any $\gamma\in\Gamma^*$, $\mu\in\Gamma_s$, $u\in{s\Gamma}$ and $\sigma\in{\rm Aut}(\GWLL)$, we can suppose
\begin{eqnarray}\label{201303050600}
\begin{array}{llll}
&&\sigma(L_0)=a_0L_0+c_0\frak{c},\ \ \ \sigma(\frak{c})=m_0\frak{c},\vs{8pt}\\
&&\sigma(L_\gamma)=\mbox{$\sum\limits_{\alpha\in\Gamma}$}a_{\gamma,\alpha}L_\alpha+
\mbox{$\sum\limits_{\nu\in\Gamma_s}$}b_{\gamma,\nu}T_\nu+c_\gamma\frak{c},\vs{8pt}\\
&&\sigma(T_\mu)=\mbox{$\sum\limits_{\alpha\in\Gamma}$}d_{\mu,\alpha}L_\alpha+
\mbox{$\sum\limits_{\nu\in\Gamma_s}$}e_{\mu,\nu}T_\nu+f_\mu\frak{c},\vs{8pt}\\
&&\sigma(G_u)=\mbox{$\sum\limits_{v\in{s\Gamma}}$}n_{u,v}G_v,
\end{array}
\end{eqnarray}
where $a_0,\,m_0,\,a_{\gamma,\alpha},\,b_{\gamma,\nu},\,c_\gamma,
\,d_{\mu,\alpha},\,e_{\mu,\nu},\,f_\mu,\,n_{u,v}\in\F$. Nearly repeating the proving process of Theorem \ref{201302281330}, we can deduce
\begin{eqnarray*}
&&a_0\neq0,\ \ \ (a_0v-u)n_{u,v}=0,\\
&&(a_0\alpha-\gamma)a_{\gamma,\alpha}=0,\ \ \ (a_0\nu-\gamma)b_{\gamma,\nu}=0,\ \ \ \gamma c_\gamma=0,\\
&&(a_0\alpha-\mu)d_{\mu,\alpha}=0,\ \ \ (a_0\nu-\mu)e_{\mu,\nu}=0,
\ \ \ \mu f_\mu=0,
\end{eqnarray*}
which imply
\begin{eqnarray*}
&&a_{\gamma,\alpha}=b_{\gamma,\nu}=c_\gamma=d_{\mu,\beta}=e_{\mu,\nu'}=f_\mu=n_{u,v}=0,
\end{eqnarray*}
for any $\gamma\in\Gamma^*$, $\beta\in\Gamma$, $\alpha\neq\varepsilon\gamma$,
$\mu,\,\nu\in\Gamma_s$, $\nu'\neq\varepsilon\mu$, $v\neq\varepsilon u$ and $\varepsilon\in\{\pm1\}$.
Then we can rewrite the identities referred in \eqref{201303050600} as follows:
\begin{eqnarray}\label{201303050601}
&&\sigma(L_\gamma)=a_{\gamma}L_{\varepsilon\gamma},\ \ \
\sigma(T_\mu)=e_{\mu}T_{\varepsilon\mu},\ \ \
\sigma(G_u)=n_{u}G_{\varepsilon u},
\end{eqnarray}
for any $\gamma\in\Gamma^*$, $\mu\in\Gamma_s$, $u\in{s\Gamma}$ with $\varepsilon\in\{\pm1\}$.

For the case $\varepsilon=1$, using \eqref{201303050601}, we can write
\begin{eqnarray}\label{201303050602}
\begin{array}{lll}
&&\sigma(\frak{c})=m_0\frak{c},\ \ \
\sigma(L_\gamma)=a_\gamma L_\gamma+c_0\delta_{\gamma,0}\frak{c},\\[8pt]
&&\sigma(T_\mu)=e_\mu T_\mu,\ \ \ \sigma(G_u)=n_uG_u,
\end{array}
\end{eqnarray}
for any $\gamma\in\Gamma$, $\mu\in\Gamma_s$, $u\in{s\Gamma}$ and the coefficients are all in $\F$.

Nearly repeating the corresponding proving process given in Theorem \ref{201302281330}, we can deduce
\begin{eqnarray*}
&&e_\mu e_{-\mu}=1,\ \ \ e_{\mu+\gamma}=e_\mu a_\gamma,\\
&&(\gamma-\alpha)(a_{\gamma+\alpha}-a_\gamma a_\alpha)=0,\\
&&\frac{\gamma^3-\gamma}{12}a_\gamma a_{-\gamma}=2\gamma c_0+\frac{\gamma^3-\gamma}{12}m_0,\\
&&\left\{\begin{array}{lll}
&\!\!\!\!\!\!n_{u}n_{v}=a_{u+v}&\mbox{if}\ \,\,u+v\in\Gamma,\\[8pt]
&\!\!\!\!\!\!n_{u}n_{v}=e_{u+v}&\mbox{if}\ \,\,u+v\in\Gamma_s,
\end{array}\right.
\end{eqnarray*}
which imply
\begin{eqnarray*}
&&c_0=0,\ \ \ m_0=1,\ \ \ a_{\gamma+\alpha}=a_\gamma a_\alpha,\\
&&e^2_\mu=a_{2\mu},\ \ \ e_{\mu+\gamma}=e_\mu a_\gamma,\ \ \
n_{u}=\left\{\begin{array}{lll}
&\!\!\!\!\!\!\varepsilon a_{u}&\mbox{if}\ \,\,u\in\Gamma,\\[8pt]
&\!\!\!\!\!\!\varepsilon e_{u}&\mbox{if}\ \,\,u\in\Gamma_s,
\end{array}\right.
\end{eqnarray*}
for all $\alpha,\,\gamma\in\Gamma$, $\mu\in\Gamma_s$ and $u\in{s\Gamma}$.
Then the identities given in \eqref{201303050602} can be rewritten as follows:
\begin{eqnarray*}
\begin{array}{lll}
&&\sigma(\frak{c})=\frak{c},\ \,\sigma(L_\gamma)=a_\gamma L_\gamma,
\ \,\sigma(T_\mu)=e_\mu T_\mu,\\[8pt]
&&\sigma(G_u)=\left\{\begin{array}{lll}
&\!\!\!\!\!\!\varepsilon a_{u}G_u&\mbox{if}\ \,\,u\in\Gamma,\\[8pt]
&\!\!\!\!\!\!\varepsilon e_{u}G_u&\mbox{if}\ \,\,u\in\Gamma_s,
\end{array}\right.
\end{array}
\end{eqnarray*}
for any $\gamma\in\Gamma$, $\mu\in\Gamma_s$, $u\in{s\Gamma}$ and $\varepsilon^2=1$,
$a_{\gamma+\alpha}=a_\gamma a_\alpha$, $e^2_\mu=a_{2\mu}$, $e_{\mu+\gamma}=e_\mu a_\gamma$.

For the case $\varepsilon=-1$, the identities given in \eqref{201303050600} and \eqref{201303050602} can be rewritten as follows:
\begin{eqnarray*}
\begin{array}{lll}
&&\sigma(\frak{c})=-\frak{c},\ \,\sigma(L_\gamma)=a_\gamma L_{-\gamma},
\ \,\sigma(T_\mu)=e_\mu T_{-\mu},\\[8pt]
&&\sigma(G_u)=\left\{\begin{array}{lll}
&\!\!\!\!\!\!-\omega a_{u}G_{-u}&\mbox{if}\ \,\,u\in\Gamma,\\[8pt]
&\!\!\!\!\!\!\omega e_{u}G_{-u}&\mbox{if}\ \,\,u\in\Gamma_s,
\end{array}\right.
\end{array}
\end{eqnarray*}
for any $\gamma\in\Gamma$, $\mu\in\Gamma_s$, $u\in{s\Gamma}$ and $\omega^2=-1$,
$a_{\gamma+\alpha}=-a_\gamma a_\alpha$, $e^2_\mu=-a_{2\mu}$, $e_{\mu+\gamma}=-e_\mu a_\gamma$.

For convenience, we introduce the following notation:
\begin{eqnarray*}
&&\omega^u_\varepsilon=\left\{\begin{array}{lll}
&\!\!\!\!\!\!\varepsilon&\mbox{if}\ \,\,\varepsilon=1,\,u\in s\Gamma,\\[8pt]
&\!\!\!\!\!\!-\omega&\mbox{if}\ \,\,\varepsilon=-1,\,u\in\Gamma,\\[8pt]
&\!\!\!\!\!\!\omega&\mbox{if}\ \,\,\varepsilon=-1,\,u\in\Gamma_s,
\end{array}\right.
\end{eqnarray*}
where $\varepsilon^2=1$ and $\omega^2=-1$.

Then for any $\gamma\in\Gamma^*$, $\mu\in\Gamma_s$, $u\in{s\Gamma}$ and $\sigma\in{\rm Aut}(\GWLL)$, the identities given in \eqref{201303050600} can be rewritten as follows:
\begin{eqnarray}\label{201303050900}
\begin{array}{lll}
&&\sigma(\frak{c})=\varepsilon\frak{c},\ \,\sigma(L_\gamma)=a_\gamma L_{\varepsilon\gamma},
\ \,\sigma(T_\mu)=e_\mu T_{\varepsilon\mu},\\[8pt]
&&\sigma(G_u)=\left\{\begin{array}{lll}
&\!\!\!\!\!\!\omega^u_\varepsilon a_{u}G_{\varepsilon u}&\mbox{if}\ \,\,u\in\Gamma,\\[8pt]
&\!\!\!\!\!\!\omega^u_\varepsilon e_{u}G_{\varepsilon u}&\mbox{if}\ \,\,u\in\Gamma_s,
\end{array}\right.
\end{array}
\end{eqnarray}
for any $\gamma\in\Gamma$, $\mu\in\Gamma_s$, $u\in{s\Gamma}$ and
\begin{eqnarray}\label{201303050901}
&&a_{\gamma+\alpha}=\varepsilon a_\gamma a_\alpha,\ \ \ e^2_\mu=\varepsilon a_{2\mu},\ \ \ e_{\mu+\gamma}=\varepsilon e_\mu a_\gamma.
\end{eqnarray}
Then this theorem follows.\BSQ

The following corollary follows immediately from Theorem \ref{201303050500} and $\Inn(\GWLL)\cong\F^*$.
\begin{coro}
${\rm Aut}(\GWLL)=\F^*\rtimes\frak{G}$.
\end{coro}

\end{document}